\documentclass[a4paper]{article}
\usepackage{graphicx, verbatim, enumerate, amsmath, amsthm, amssymb, amstext, color, amsfonts}
\usepackage[format=hang, font=small, labelfont=it, margin=15mm]{caption}
\usepackage[arrow, matrix, curve]{xy}
\usepackage[applemac]{inputenc}

\theoremstyle{plain}
\newtheorem{theo}{Theorem}

\newtheorem{lem}[theo]{Lemma}

\newtheorem{cor}[theo]{Corollary}

\theoremstyle{definition}

\newtheorem{example}{Example}

\newcommand{\N}{\ensuremath{\mathbb{N}}}
\newcommand{\Z}{\ensuremath{\mathbb{Z}}}

\newcommand{\sm}{\ensuremath{\smallsetminus}}

\newcommand{\sub}{\subseteq}

\newcommand{\Gtilde}{{\widetilde G}}
\newcommand{\GPtilde}{{\widetilde G_\Psi}} 
\newcommand{\GOtilde}{{\widetilde G_\Omega}} 

\newcommand{\tst}{topological spanning tree}

\newcommand{\COMMENT}[1]{}

\newenvironment{txteq}
  {
    \begin{equation}
    \begin{minipage}[c]{0.85\textwidth} 
    \em                                
  }
  {\end{minipage}\end{equation}\ignorespacesafterend}
\newenvironment{txteq*}
  {
    \begin{equation*}
    \begin{minipage}[c]{0.85\textwidth} 
    \em                                
  }
  {\end{minipage}\end{equation*}\ignorespacesafterend}

\DeclareMathOperator{\im}{im}

\advance\textheight by -1cm

\author{Reinhard Diestel \and Julian Pott}
\title{Dual trees must share their ends}
\date{}

\begin{document}
\maketitle

\begin{abstract}
\noindent
We extend to infinite graphs the matroidal characterization of finite graph duality, that two graphs are dual iff they have complementary spanning trees in some common edge set. The naive infinite analogue of this fails.

The key in an infinite setting is that dual trees must share between them not only the edges of their host graphs but also their ends: the statement that a set of edges is acyclic and connects all the vertices in one of the graphs iff the remaining edges do the same in its dual will hold only once each of the two graphs' common ends has been assigned to one graph but not the other, and `cycle' and `connected' are interpreted topologically in the space containing the respective edges and precisely the ends thus assigned.

This property characterizes graph duality: if, conversely, the spanning trees of two infinite graphs are complementary in this end-sharing way,%
   \COMMENT{}
   the graphs form a dual pair.
\end{abstract}

\section{Introduction}
It is well known (and not hard to see) that two finite graphs are dual if and only if they can be drawn with a common abstract set of edges so that the edge sets of the spanning trees of one are the complements of the edge sets of the spanning trees of the other:

\begin{theo}\label{finite}
Let $G=(V,E)$ and $G^*=(V^*,E)$ be connected finite graphs with the same abstract edge set. Then the following statements are equivalent:
\begin{enumerate}[\rm (i)]
\item $G$ and $G^*$ are duals of each other.
\item Given any set $F\sub E$, the graph $(V,F)$ is a tree if and only if $(V^*,F^\complement)$ is a tree.
\end{enumerate}
\end{theo}

For infinite dual graphs $G$ and~$G^*$ (see~\cite{duality}), Theorem~\ref{finite}\,(ii) will usually fail: when $(V,F)$ is a spanning tree of~$G$, the subgraph $(V^*,F^\complement)$ of $G^*$ will be acyclic but may be disconnected. For example, consider as $G$ the infinite $\Z\times \Z$ grid, and let $F$ be the edge set of any spanning tree containing a two-way infinite path, a {\em double ray}~$R$. Then the edges of $R$ will form a cut in~$G^*$, so $(V^*,F^\complement)$ will be disconnected.

Although the graphs $(V^*,F^\complement)$ in this example will always be disconnected, they become arc-connected (but remain {\em acirclic\/}) when we consider them as closed subspaces of the topological space obtained from $G^*$ by adding its end. Such subspaces are called {\em\tst s\/}; they provide the `correct' analogues in infinite graphs of spanning trees in finite graphs for numerous problems, and have been studied extensively~\cite{RDsBanffSurvey, TST}. For $G = \Z\times\Z$, then, the complements of the edge sets of ordinary spanning trees of $G$ form \tst s in~$G^*$, and vice versa (as $\Z\times\Z$ is self-dual). 

It was shown recently in the context of infinite matroids~\cite{BruhnDiestelMatroidsGraphs} that this curious phenomenon is not specific to this example but occurs for all dual pairs of graphs: neither ordinary nor topological spanning trees permit, by themselves, an extension of Theorem~\ref{finite} to infinite graphs, but as soon as one notion is used for $G$ and the other for~$G^*$, the theorem does extend. The purpose of this paper is to explain this seemingly odd phenomenon by a more general duality for graphs with ends, in which it appears as merely a pair of extreme cases.

It was shown in~\cite{BruhnSteinEndDuality} that 2-connected%
   \COMMENT{}
   dual graphs do not only have the `same' edges but also the `same' ends: there is a bijection between their ends that commutes with the bijection between their edges so as to preserve convergence of edges to ends. Now if $G$ and $G^*$ are dual 2-connected graphs with edge sets~$E$ and end sets~$\Omega$, our result is that if we specify {\em any\/} subset $\Psi$ of~$\Omega$ and consider topological spanning trees of~$G$ in the space obtained from $G$ by adding only the ends in~$\Psi$, then Theorem~\ref{finite}\,(ii) will hold if the subgraphs $(V^*,F^\complement)$ of~$G^*$ are furnished with precisely the ends in~$\Omega\setminus\Psi$. (Our earlier example is the special case of this result with either $\Psi=\emptyset$ or $\Psi=\Omega$.) And conversely, if the spanning trees of two graphs $G$ and $G^*$ with common edge and end sets complement each other in this way for some---equivalently, for every---subset $\Psi$ of their ends then $G$ and $G^*$ form a dual pair.

Here, then, is the formal statement of our theorem. A graph~$G$ is \emph{finitely separable} if any two vertices can be separated by finitely many edges; as noted by Thomassen~\cite{thomassen80, thomassen82}, this slight weakening of local finiteness is necessary for any kind of graph duality to be possible. The \emph{$\Psi$-trees\/} in~$G$, for subsets $\Psi$ of its ends, will be defined in Section~\ref{definitions}. Informally, they are the subgraphs that induce no cycle or topological circle in the space which $G$ forms with the ends in~$\Psi$ (but no other ends) and connect any two vertices by an arc in this space.

\begin{theo}\label{main}
Let $G = (V, E, \Omega)$ and $G^* = (V^*\!, E, \Omega)$ be finitely separable $2$-connected graphs with the same edge set $E$ and the same end set~$\Omega$, in the sense of~{\rm \cite{BruhnSteinEndDuality}}. Then the following assertions are equivalent:
\begin{enumerate}[\rm (i)]
\item $G$ and $G^*$ are duals of each other.%
   \COMMENT{}
\item For all $\Psi\sub\Omega$ and $F\sub E$ the following holds: $F$ is the edge set of a $\Psi$-tree in~$G$ if and only if $F^\complement$ is the edge set of a $\Psi^\complement$-tree in~$G^*$.
\item There exists a set $\Psi\sub\Omega$ such that for every $F\sub E$ the following holds: $F$ is the edge set of a $\Psi$-tree in~$G$ if and only if $F^\complement$ is the edge set of a $\Psi^\complement$-tree in~$G^*$.
\end{enumerate}
\end{theo}

\goodbreak

Setting $\Psi=\emptyset$ in (ii) and~(iii) as needed,
   \COMMENT{}
   we reobtain the following result from~\cite{BruhnDiestelMatroidsGraphs}:

\begin{cor}\label{maincor}
Two 2-connected and finitely separable graphs $G=(V,E,\Omega)$ and $G^*=(V^*,E,\Omega)$%
    \COMMENT{}
   are dual if and only if the following assertions are equivalent for every $F\sub E$:
\begin{enumerate}[\rm (i)]
\item $F$ is the edge set of a spanning tree of $G$;
\item $F^\complement$ is the edge set of a topological spanning tree of~${G^*}$.\qed
\end{enumerate}
\end{cor}

We shall prove Theorem~\ref{main}, extended by another pair of equivalent conditions in terms of circuits and bonds, in Sections \ref{basics}--\ref{proofofmain}.

\section{Definitions and basic facts}\label{definitions}

All the graphs we consider in this paper will be \emph{finitely separable}, that is, any two vertices can be separated by finitely many edges.%
   \COMMENT{}

We think of a \emph{graph} as a triple $(V,E,\Omega)$ of disjoint sets, of \emph{vertices}, \emph{edges}, and \emph{ends}, together with a map $ E\to V\cup [V]^2$ assigning to every edge either one or two vertices, its \emph{endvertices}, and another map mapping the ends bijectively to the equivalence classes of \emph{rays} in the graph, its 1-way infinite paths, where two rays are \emph{equivalent} if they cannot be separated by finitely many vertices. In particular, our `graphs' may have multiple edges and loops. For the complement of $F$ in~$E$, and of~$\Psi$ in~$\Omega$, we write $F^\complement$ and~$\Psi^\complement$, respectively.

Let $G = (V,E,\Omega)$ be a graph, and let $X$ be the topological%
   \COMMENT{}
   $1$-complex formed by its vertices and edges. In~$X$, every edge is a topological copy of~$[0,1]$ inheriting also its metric. We denote the topological interior of an edge $e$ by~$\mathring e$, and for a set $F\sub E$ of edges we write $\mathring F:= \bigcup_{e\in F}\mathring e$.

Let us define a new%
   \COMMENT{}
   topology on~$X \cup \Omega$, to be called \emph{\textsc{VTop}}. We do this by specifying a neighbourhood basis for every point.
For points $x\in X$ we declare as open the open $\epsilon$-balls around $x$ in~$X$ with $0<\epsilon<\delta$, where $\delta$ is the distance from $x$ to a closest vertex $v\ne x$. For points $\omega\in\Omega$, note that for every finite set $S\sub V$ there is a unique component $C = C(S,\omega)$ of $G-S$ that contains a ray from~$\omega$.
Let $\hat C = \hat C(S,\omega)\sub X\cup\Omega$ be the set of all the vertices and inner points of edges contained in or incident with~$C$, and of all the ends represented by a ray in~$C$. We declare all these sets $\hat C$ as open, thus obtaining for $\omega$ the neighbourhood basis
 $$\big\{\hat C(S,\omega)\sub X\cup \Omega : S\sub V,\ |S|<\infty\big\}.$$
We write $|G|$ for the topological space on $X\cup\Omega$ endowed with this topology.%
   \footnote{This differs a little from the definition of $|G|$ in~\cite{DiestelBook10noEE} when $G$ is not locally finite.}
In topological contexts we shall also write $G$ for the subspace $|G|\sm\Omega$. (This  has the same points as~$X$, but a different topology unless $G$ is locally finite.)%
  \COMMENT{}

If $\omega$ and $S$ are as above, we say that $S$ {\em separates\/} $\omega$ in $G$ from all the ends that have no ray in $C(S,\omega)$ and from all vertices in $G - C(S,\omega) - S$.

A vertex $v$ \emph{dominates} an end $\omega$ if $G$ contains infinitely many paths from $v$ to some ray in $\omega$ that pairwise meet only in~$v$. When this is the case we call $v$ and~$\omega$ {\em equivalent\/}; let us write $\sim$ for the equivalence relation on $V\cup\Omega$ which this generates. Note that since $G$ is finitely separable, no two vertices will be equivalent under~$\sim\,$: every non-singleton equivalence class consists of one vertex and all the ends it dominates. A~vertex and an end it dominates have no disjoint neighbourhoods in~$|G|$. But two ends always have disjoint neighbourhoods, even if they are dominated by the same vertex.%
   \COMMENT{}

For sets $\Psi\sub\Omega$ of ends, we shall often consider the subspace
 $$|G|_\Psi:=|G|\sm\Psi^\complement$$
and its quotient space
 $${\GPtilde}:=|G|_\Psi/\!\!\sim\,,$$%
   \COMMENT{}
   whose topology we denote by \emph{$\Psi$-\textsc{Top}}. For $\Psi=\Omega$ we obtain an identification space
 $$\widetilde G:={\GOtilde}$$
 that readers may have met before; its topology is commonly denoted as~\textsc{ITop}. We usually write $[x]_\Psi$ for the equivalence class of $x$ in~$|G|_\Psi$, and $[x]$ for its class in~$\Gtilde$.

As different vertices are never equivalent, the vertices of $G$ determine distinct $\sim$-classes, which we call the {\em vertices\/} of~${\GPtilde}$. All other points of $\GPtilde$ are singleton classes $\{x\}$, with $x$ either an inner point of an edge or an undominated end in~$\Psi$. We will not always distinguish $\{x\}$ from $x$ in these cases, i.e., call these~$x$ also {\em inner point of edges\/} or {\em ends\/} of~$\GPtilde$.%
   \COMMENT{}

Note that if $\Psi$ contains a dominated end then $|G|_\Psi$ will fail to be Hausdorff, and if $\Psi^\complement\ne\emptyset$ then ${\GPtilde}$ will fail to be compact. But we shall see that $\GPtilde$ is always Hausdorff (Corollary~\ref{HD}), and if $G$ is 2-connected then $\widetilde G$ is compact~\cite{diestelESST}.

Rather than thinking of $\GPtilde$ as a quotient space as formally defined above, we may think of it informally as formed from the topological space~$G$ in three steps:
   \begin{itemize}
   \item add the undominated ends from~$\Psi$ as new points, and make their rays converge to them;
   \item make the rays from any dominated end in~$\Psi$ converge to their unique dominating vertex;
   \item let the rays of ends in $\Psi^\complement$ go to infinity without converging to any point.
   \end{itemize}

The diagram in Figure~\ref{diagram} shows the relationship between the spaces just defined.
The subspace inclusion $\iota\colon |G|_\Psi\to |G|$ and the quotient projections $\pi\colon |G|\to \widetilde G$ and $\pi_\Psi\colon |G|_\Psi\to{\GPtilde}$ are canonical, and ${\sigma_\Psi}\colon {\GPtilde}\to\widetilde G$ is defined so as to make the diagram commute: it sends an equivalence class $[x]_\Psi\in {\GPtilde}$ to the class $[x]\in \widetilde G$ containing it.

\begin{figure}[htbp]
\begin{center}
\mbox{
\begin{xy}
  \xymatrix{
      |G|_\Psi \ar[r]^\iota \ar[d]_{\pi_\Psi}    &   |G| \ar[d]^\pi  \\
      {\GPtilde} \ar[r]_{\sigma_\Psi}             &   \widetilde G  
  }
\end{xy}
}
\end{center}
\caption{Spaces with ends, and their quotient spaces}\label{diagram}\vskip-12pt\vskip0pt
\end{figure}

\goodbreak
   
Since $G$ is finitely separable and hence no end is dominated by more than one vertex, ${\sigma_\Psi}$~is injective: ${\sigma_\Psi}([x]_\Psi) = [x]\in\widetilde G$ is obtained from $[x]_\Psi$ simply by adding those ends of $\Psi^\complement$ that are dominated by a vertex in~$[x]_\Psi$. As $|G|_\Psi$ carries the subspace topology induced from~$|G|$, it is also easy to check that ${\sigma_\Psi}$ is continuous. Its inverse~$\sigma_\Psi^{-1}$ can fail to be continuous; see Example~\ref{psicontinuity} below.

The subtle differences between $|G|_\Psi$ and $\GPtilde$ will often be crucial in this paper. But when they are not, we may suppress them for simplicity of notation. For example, given a subgraph $H$ of $G$ we shall speak of the {\em closure of $H$ in~$\GPtilde$\/} and mean the obvious thing:  the closure in $\GPtilde$ of its subspace $\pi_\Psi (H')$, where $H'$ is $H$ viewed as a subspace of~$|G|_\Psi\sub |G|$.

By a {\em circle\/} in a topological space $X$ we mean a topological embedding ${S^1\to X}$, or its image. Since circles are compact and $\widetilde G$ is Hausdorff,%
   \COMMENT{}
   ${\sigma_\Psi}$ maps circles in $\GPtilde$ to circles in~$\Gtilde$. Conversely, circles in $\Gtilde$ that use only ends in~$\Psi$ define circles in~$\GPtilde$; this will be shown in Lemma~\ref{circles}. The set of all the edges contained in a given circle in~$\Gtilde_\Psi$ will be called a $\Psi$-{\em circuit\/} of~$G$; for $\Psi=\Omega$ we just speak of {\em circuits} of~$G$. We shall not consider `circuits' of circles in $|G|$ or~$|G|_\Psi$.

As with circles, we use the term {\em path\/} in topological contexts both for continuous maps from $[0,1]$, not necessarily injective, and for their images. For example, if $A$ and $B$ are the images of paths $\varphi,\varphi'\colon[0,1]\to \widetilde G$ with endpoints $x = \varphi(0)$ and $y = \varphi(1) = \varphi'(0)$ and $z = \varphi'(1)$, we write $xAyBz$ for the `$x$--$y$ path' in $\Gtilde$ that is the image of the concatenation of the paths $\varphi$ and~$\varphi'$. Note that, since $\GPtilde$ is Hausdorff, every path in $\GPtilde$ between two points $x$ and~$y$ contains an $x$--$y$ arc \cite[p.\,208]{ElemTop}.

A subspace of ${\GPtilde}$ that is the closure in $\GPtilde$ of the union of all the edges it contains is a \emph{standard subspace of~${\GPtilde}$}. Circles in $\Gtilde_\Psi$ are examples of standard subspaces; this was shown in~\cite{TST} for~$\Gtilde$, and follows for arbitrary~$\Psi$ from Lemma~\ref{finsepext} below.%
   \COMMENT{}
   A~standard subspace of $\GPtilde$ that contains no circle is a  \emph{$\Psi$-forest} of~$G$. A~\hbox{$\Psi$-forest} is {\em spanning\/} if it contains all the vertices of~$\GPtilde$. Note that, being closed, it then also contains all the ends of~$\GPtilde$. A spanning arc-connected $\Psi$-forest of $G$ is a \emph{$\Psi$-tree} of~$G$.

Thus, the $\emptyset$-trees of~$G$ are precisely its (ordinary) spanning trees, while its $\Omega$-trees are its \emph{\tst s}, the arc-connected standard subspaces of~$\Gtilde$ that contain all the vertices of $G$ but no topological circle.

\begin{figure}[htpb]
\centering\vskip-12pt
   	  \includegraphics{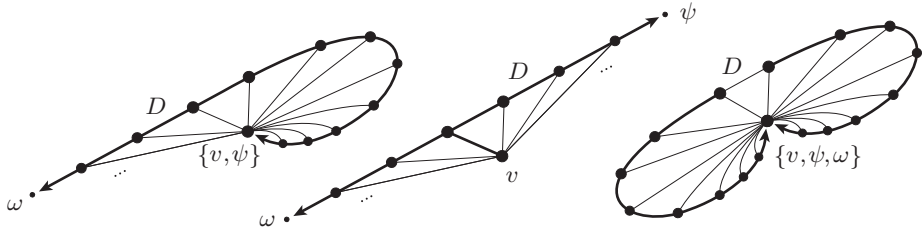}
   	  \caption{$\Psi$-trees for $\Psi=\{\psi\}$, \  $\Psi=\emptyset$ and $\Psi = \{\omega,\psi\}$}
   \label{Psitrees}\vskip-6pt\vskip0pt
   \end{figure}

\begin{example}
Let $G $ be obtained from a double ray~$D$ by adding a vertex~$v$ adjacent to all of~$D$. This graph $G$ has two ends, $\omega$ and~$\psi$ say, both dominated by~$v$. The closure in~$\Gtilde$ of the edges of $D$ is a circle containing the `vertex' $[v] = \{v,\omega,\psi\}$ of~$\Gtilde$, even though $v$ does not lie on~$D$. However for $\Psi = \{\psi\}$ the closure in~$\Gtilde_\Psi$ of the same set of edges is not a circle but homeomorphic to a half-open interval. It thus is a $\Psi$-tree, and even a spanning one, since $v$ and $\psi$ are both elements of its `vertex' $\{v,\psi\}$ and it also contains all the other vertices of~$G$. The closure of the edges of $D$ in $\Gtilde_\emptyset$, on the other hand, is a $\emptyset$-tree but not a spanning one, since $v$ lies in none of its points. Figure~\ref{Psitrees} shows a $\Psi$-tree for each choice of $\Psi$ in this example.
\end{example}

If $G$ and $G^*$ are graphs with the same edge set, and such that the bonds of $G^*$ are precisely the circuits of~$G$, then $G^*$ is called a {\em dual\/} of~$G$. If the {\em finite\/} bonds of $G^*$ are precisely the finite circuits of~$G$, then $G^*$ is a {\em finitary dual\/} of~$G$. Clearly, duals are always finitary duals. For finitely separable graphs, as considered here, the converse is also true \cite[Lemmas 4.7--4.9]{duality}. If $G^*$ is a dual of~$G$, then $G$ is a dual of~$G^*$~\cite[Theorem~3.4]{duality}. Finally, $G$~has a dual if and only if it is planar~\cite{duality}.%
   \COMMENT{}

\section{Lemmas}\label{basics}

Our main aim in this section is to prove some fundamental lemmas about the spaces $|G|$, $|G|_\Psi$, $\Gtilde$ and~$\Gtilde_\Psi$ defined in Section~\ref{definitions}: about their topological properties, and about their relationship to each other. Throughout the section, let $G=(V,E,\Omega)$ be a fixed finitely separable graph, and $\Psi\sub \Omega$ a fixed set of ends.

Before we get to these topological fundamentals, let us show that $\Psi$-trees always exist, and prove an easy lemma about how they relate to finite circuits and bonds. As to the existence of $\Psi$-trees, we can even show that there are always rather special ones: $\Psi$-trees that are connected not only topologically through their ends, but also as graphs:

\begin{lem}\label{existence}
If $G$ is connected, it has a spanning tree $T$ whose closure in $\GPtilde$ is a $\Psi$-tree.
\end{lem}

\begin{proof}
It was shown in~\cite[Thm.~6.3]{duality} that $G$ has a spanning tree $T $ whose closure $\overline T$ in $\Gtilde$ contains no circle. Let ${\overline T_\Psi}$ denote the closure of~$T$ in~$\GPtilde$. Then $\overline T = \sigma_\Psi (\overline T_\Psi)$. Since circles in $\Gtilde_\Psi$ define circles in~$\Gtilde$ (by composition with~$\sigma_\Psi$), ${\overline T_\Psi}$~contains no circle either.

\goodbreak

For a proof that ${\overline T_\Psi}$ is arc-connected it suffices to show that every undominated end $\psi\in\Psi$ contains a ray $R\sub T$: then the arc $\pi_\Psi(T)\sub \overline T_\Psi$ connects the end $\{\psi\}\in \overline T_\Psi$ to a vertex, while all the vertices of $\overline T_\Psi$ are connected by~$T$. Pick a ray $R'\in\psi$ in~$G$, say $R' = v_0 v_1\dots$. By the star-comb lemma \cite[Lemma 8.2.2]{DiestelBook10noEE}, the connected graph $\bigcup_{n\in\N} v_n T v_{n+1}$ contains a subdivided infinite star with leaves in~$R'$ or an infinite comb with teeth in~$R'$. As $\psi$ is not dominated,%
   \COMMENT{}
   we must have a comb. The back $R\sub T$ of this comb is a ray equivalent to~$R'$ that hence lies in~$\psi$.

Being acirclic, arc-connected and spanning, ${\overline T_\Psi}$~is a $\Psi$-tree.
\end{proof}

\begin{lem}\label{matroidlemma}\COMMENT{}
Assume that $G$ is connected, and let $F\sub E$ be a finite set of edges.
\begin{enumerate}[\rm (i)]
\item $F$ is a circuit if and only if it is not contained in the edge set of any $\Psi$-tree and is minimal with this property.
\item $F$ is a bond if and only if it meets the edge set of every $\Psi$-tree and is minimal with this property.
\end{enumerate}
\end{lem}

\begin{proof}
(i) Assume first that $F$ is a circuit. Then $F$ is not contained in any $\Psi$-tree; let us show that every proper subset of $F$ is. We do this by showing the following more general fact:
\begin{txteq}\label{extendarc}
Every finite set $F'$ of edges not containing a circuit extends to a spanning tree of $G$ whose closure in $\GPtilde$ is a $\Psi$-tree.
\end{txteq}
To prove~\eqref{extendarc}, consider a spanning tree $T$ of~$G$ whose closure in $\GPtilde$ is a $\Psi$-tree (Lemma~\ref{existence}). Choose it with as many edges in $F'$ as possible. Suppose it fails to contain an edge $f\in F'$. Adding $f$ to $T$ creates a cycle $C$ in $T+ f$, which by assumption also contains an edge $e\notin F'$. As $C$ is finite, it is easy to check that $T+f-e$ is another spanning tree whose closure is a $\Psi$-tree.%
   \COMMENT{}
   This contradicts our choice of~$T$.

Conversely, if $F$ is not contained in any $\Psi$-tree, then by~\eqref{extendarc} it contains a circuit. If, in addition, it is minimal with the first property, it will in fact be that circuit, since we could delete any other edge without making it extendable to a $\Psi$-tree.

(ii) If $F$ is a cut, $F=E(V_1,V_2)$ say, then the closures of $G[V_1]$ and $G[V_2]$ in $\GPtilde$ are disjoint open subsets of $\GPtilde\sm\mathring F$, so this subspace cannot contain a $\Psi$-tree. Thus, $F$~meets the edge set of every $\Psi$-tree.

If $F$ is even a bond, then both $V_1$ and $V_2$ induce connected subgraphs. By Lemma~\ref{existence}, these have spanning trees $T_i$ ($i=1,2$) whose closures in~$\GPtilde$%
   \COMMENT{}
   are arc-connected and contain no circle.%
   \footnote{We are applying Lemma~\ref{existence} in the subgraphs~$G[V_i]$. But since $F$ is finite, the spaces $\widetilde{G[V_i]}_{\Psi_i}$ are canonically embedded in~$\GPtilde$.}
   For every edge $f\in F$, the closure $\overline T_\Psi$ of $T:= (T_1\cup T_2) + f$ in $\GPtilde$ then is a $\Psi$-tree of~$G$: it still contains no circle, because no arc in $\overline T_\Psi\sm\mathring f$ can cross the finite cut $F$ from which it contains no edge (as above).%
   \COMMENT{}
   So $F$ is minimal with the property of meeting the edge set of every $\Psi$-tree.

Conversely, let us assume that $F$ meets the edge set of every $\Psi$-tree, and show that $F$ contains a bond.%
   \COMMENT{}
   Let $T$ be a spanning tree of $G$ whose closure in $\GPtilde$ is a $\Psi$-tree (Lemma~\ref{existence}), chosen with as few edges in $F$ as possible. By assumption, $T$~has an edge $f$ in~$F$. If the bond $B$ of $G$ between the two components of $T-f$ contains an edge $e\notin F$, then $T-f+e$ is another spanning tree whose closure is a $\Psi$-tree (as before) that contradicts our choice of~$T$. So $B$ contains no such edge~$e$ but is contained in~$F$.

If $F$ is minimial with the property of containing an edge from every $\Psi$-tree, it must be equal to the bond it contains. For by the forward implication of~(ii) already proved, any other edge could be deleted from $F$ without spoiling its property of meeting the edge set of every $\Psi$-tree.
\end{proof}

We begin our study of the spaces introduced in Section~\ref{definitions} by showing that finite separability extends from $G$ to~${\GPtilde}$:

\begin{lem}\label{finsepext}
For every two points $p,q\in{\GPtilde}$ that are not inner points of edges there exists a finite set $F$ of edges such that $p$ and $q$ lie in disjoint open sets of~$\GPtilde\sm \mathring F$ whose union is $\GPtilde\sm \mathring F$.
\end{lem}

\begin{proof}
Let us write $p=[x]_\Psi$ and~$q=[y]_\Psi$, where $x$ and $y$ are either vertices or undominated ends of~$G$. We shall find a finite cut $F$ of $G$, with bipartition $(X,Y)$ of $V$ say, such that $x\in\overline X$ and $y\in\overline Y$, where $\overline X$ and $\overline Y$ denote closures of $X$ and $Y$ in~$|G|_\Psi$. Since $F$ is finite, $\overline X$ and $\overline Y$ then partition of~$|G|_\Psi\sm\mathring F$ into disjoint open sets that are closed under equivalence, so their projections under $\pi_\Psi$ partition $\GPtilde\sm\mathring F$ into disjoint open sets containing $p$ and~$q$, respectively.%
	\COMMENT{}

If $x$ and $y$ are vertices, then $F$ exists by our assumption that $G$ is finitely separable.%
	\COMMENT{}
Suppose now that $y$ is an end. Let us find a finite set $S\not\owns x$ of vertices that separates $x$ from $y$ in~$G$.%
   \COMMENT{}
   If $x$ is another end, then $S$ exists since $x\ne y$. If $x$ is a vertex, pick a ray $R\in y$. If there is no $S$ as desired, we can inductively find infinitely many independent $x$--$R$ paths in~$G$,%
   \COMMENT{}
   contradicting the fact that $y$ is undominated.

Having found $S$, consider the component $C := C(S,y)$ of $G-S$. For each $s\in S$ we can find a finite set $S_s\sub C$ of vertices separating $s$ from~$y$ in the subgraph of $G$ spanned by $C$ and~$s$, since otherwise $s$ would dominate~$y$ (as before). Let $S':= \bigcup_{s\in S} S_s$; this is a finite set of vertices in $C$ that separates all the vertices of $S$ from $y$ in~$G$. Since $G$ is finitely separable, there is a finite set $F$ of edges separating $S$ from $S'$ in~$G$. Choose $F$ minimal. Then, assuming without loss of generality that $G$ is connected, every component of $G-F$ meets exactly one of the sets $S$ and~$S'$. Let $X$ be the set of vertices in components meeting~$S$, and let $Y$ be the set of vertices in components meeting~$S'$. Then $(X,Y)$ is a partition of $G$ crossed by exactly the edges in~$F$, and it is easy to check that $F$ has the desired properties.%
   \COMMENT{}
   \end{proof}

It was proved in~\cite{TST}, under a weaker assumption than finite separability (just strong enough that $\Gtilde$ can be defined without identifying distinct vertices)%
   \COMMENT{}
   that $\Gtilde$ is Hausdorff. For finitely separable graphs, as considered here, the proof is much simpler and extends readily to~$\GPtilde$:

\begin{cor}\label{HD}
$\GPtilde$ is Hausdorff.
\end{cor}

\begin{proof}%
   \COMMENT{}
Finding disjoint open neigbourhoods for distinct points $p,q\in\GPtilde$ is easy if one of them is an inner point of an edge. Assume that this is not the case, let $F$, $\overline X$ and $\overline Y$ be defined as in Lemma~\ref{finsepext} and its proof, and let $S$ be the (finite) set of vertices incident with an edge in~$F$. Then $p\sub \overline X$ and $q\sub\overline Y$. Any end $\psi\in p$ has a basic open neighbourhood $\hat C(S,\psi)\sm\Psi^\complement$ in~$|G|_\Psi$ that is a subset of~$\overline X\sm S$. Write $O_p$ for the union of all these neighbourhoods, together with a small open star neighbourhood of the vertex in $p$ if it exists.%
   \COMMENT{}
   Define $O_q$ similarly for~$q\sub\overline Y$. Then $\pi_\Psi (O_p)$ and $\pi_\Psi (O_q)$ are disjoint open neighbourhoods of $p$ and~$q$ in~$\GPtilde$.
\end{proof}

Our next aim is to select from the basic open neighbourhoods ${\hat C(S,\omega)\sm\Psi^\complement}$ in $|G|_\Psi$ of ends ${\omega\in\Psi}$ some `standard' neighbourhoods that behave well under the projection~$\pi_\Psi$ and still form neighbourhood bases of these points~$\omega$. Ideally, we would like to find for every end $\omega\in\Psi$ a basis of open neighbourhoods that are closed under~$\sim\,$. That will not be possible, since ends $\omega'\ne\omega$ equivalent to~$\omega$ can be separated topologically from~$\omega$. But we shall be able to find a basis of open neighbourhoods of $\omega$ that will be closed under~$\sim$ for all points other than $\omega$ itself. Then the union of all these neighbourhoods, one for every end $\omega'\sim\omega$, plus an open star neighbourhood of their common dominating vertex, will be closed under~$\sim\,$, and will thus be the pre-image of an open neighbourhood of $\pi_\Psi (\omega) = [\omega]_\Psi$ in~$|G|_\Psi$.

Given a bond $F=E(V_1,V_2)$ of~$G$%
   \COMMENT{}
   and an end $\omega\in\Psi$ that lies in the $|G|$-closure of $V_1$ but not of~$V_2$, let
 $$\hat C_\Psi(F,\omega)\sub|G|_\Psi$$
 denote the union of the $|G|_\Psi$-closure of $G[V_1]$ with~$\mathring F$. For every vertex $v\in V_2$ we also call $F$ a \emph{$v$--$\omega$ bond\/}. Note that $\hat C_\Psi(F,\omega)$ depends only on $F$ and~$\omega$: since $F$ is a bond, $G-F$ has only two components, so $V_1$ and $V_2$ can be recovered from $F$ and~$\omega$. Note also that every ray in $\omega$ has a tail in~$\hat C_\Psi(F,\omega)$,%
   \COMMENT{}
  so if it starts at~$v$ it must have an edge in~$F$.

If $v\in V_2$ is an endvertex of all but finitely many of the edges in~$F$, we say that $F$ is \emph{$v$-cofinite\/}. Then the set $S$ of endvertices of $F$ in~$V_2$ is finite and separates $\omega$ from $V_2\sm S$.

\begin{lem}\label{standardnbhds}
Let $\omega\in\Psi$ be an end, and $v\in V$ a vertex.
\begin{enumerate}[\rm (i)]
\item If $\omega$ is undominated, then the sets $\{\,\hat C_\Psi(F,\omega)\mid F\text{ is a finite bond of }G\,\}$ form a basis of open neighbourhoods of~$\omega$ in~$|G|_\Psi$.
\item If $\omega$ is dominated by~$v$, then the sets
 $$\{\,\hat C_\Psi(F,\omega)\mid F\text{ is a $v$-cofinite $v$--$\omega$ bond}\,\}$$
form a basis of open neighbourhoods of~$\omega$ in~$|G|_\Psi$.
\end{enumerate}
\end{lem}

\begin{proof}
(i) As $F$ is finite, so is the set $S$ of its endvertices in~$V_2$. Since $F$ is a bond, $G[V_1]$~is connected. Hence $\hat C_\Psi(F,\omega)$ equals $\hat C(S,\omega)\sm\Psi^\complement$, which is a basic open neighbourhood of~$\omega$ in~$|G|_\Psi$.
Conversely, we need to find for any finite set $S\sub V$, without loss of generality connected,%
	\footnote{The sets $\hat C(S,\omega)$ with $S$ connected in~$G$ also form a neighbourhood basis of~$\omega$ in~$|G|$, since every finite set $S$ of vertices extends to a finite connected set.}
   a finite bond $F$ such that $\hat C_\Psi(F,\omega)\sub \hat C(S,\omega)$. As no vertex dominates~$\omega$, there is a finite connected set $S'$ of vertices of~$C(S,\omega)$ that separates $S$ from~$\omega$ in~$G$. (Otherwise we could inductively construct an infinite set of disjoint paths in $C(S,\omega)$ each starting at a vertex adjacent to~$S$ and ending on some fixed ray $R\in\omega$; then infinitely many of the starting vertices of these paths would share a neighbour in~$S$, which would dominate~$\omega$.) As $G$ is finitely separable, there is a finite set of edges separating $S$ from~$S'$ in~$G$. As both $S$ and $S'$ are connected, choosing this set minimal ensures that it is a bond. This bond $F$ satisfies $\hat C_\Psi(F,\omega)\sub \hat C(S,\omega)$.%
   \COMMENT{}

(ii) Although $F$ is infinite now, the set $S$ of its endvertices in~$V_2$ is finite. Hence $\hat C_\Psi(F,\omega)$ is a basic open neighbourhood of~$\omega$ in~$|G|_\Psi$,%
   \COMMENT{}
   as in the proof of~(i). Conversely, let a finite set $S\sub V$ be given; we shall find~a $v$-cofinite $v$--$\omega$ bond $F$ such that $\hat C_\Psi(F,\omega)\sub \hat C(S,\omega)$. The sets $\hat C(T,\omega)$ such that $v\in T$ and both $T-v$ and $T$ are connected in~$G$ still form a neighbourhood basis for $\omega$ in~$|G|$,%
   \COMMENT{}
   so we may assume that $S$ has these properties. As in the proof of~(i), there is a finite connected set $S'$ of vertices in~$C(S,\omega)$ that separates $S-v$ from~$\omega$ in~$G-v$, because $\omega$ is not dominated in~$G-v$. As $G-v$ is finitely separable, there is a finite bond $F= E(V_1,V_2)$ of $G-v$ that separates $S-v$ from~$S'$, with $S-v\sub V_2$ say. Then $F':=E(V_1,V_2\cup\{v\})$ is a $v$-cofinite $v$--$\omega$ bond in $G$ with $\hat C_\Psi(F',\omega)\sub \hat C(S,\omega)$, as before.
\end{proof}

Let us call the open neighbourhoods $\hat C_\Psi (F,\omega)$ from Lemma~\ref{standardnbhds} the \emph{standard neighbourhoods} in $|G|_\Psi$ of the ends $\omega\in\Psi$. For points of $|G|_\Psi$ other than ends, let their \emph{standard neighbourhoods} be their basic open neighbourhoods defined in Section~\ref{definitions}.

Trivially, standard neighbourhoods of vertices and inner points of edges are closed under~$\sim\,$. Our next lemma says that standard neighbourhoods of ends are nearly closed under~$\sim\,$, in that only the end itself may be equivalent to points outside: to a vertex dominating it, and to other ends dominated by that vertex.\looseness=-1

\begin{lem}\label{stdnhd}
If $\hat C = \hat C_\Psi(F,\omega)$ is a standard neighbourhood of $\omega\in \Psi$ in $|G|_\Psi$, then $[x]_\Psi\sub \hat C$ for every $x\in \hat C\sm[\omega]_\Psi$.
\end{lem}

\begin{proof}
Let $S$ be the finite set of vertices not in $\hat C$ that are incident with an edge in~$F$. Suppose, for a contradiction, that there are points $x\sim y$ in $|G|_\Psi$ such that $x\in \hat C\sm[\omega]$ but $y\notin \hat C\sm[\omega]$. Since the unique vertex in the $\sim_\Psi$-class of $x$ and $y$ lies either in $\hat C\sm [\omega]$ or not, we may assume that either $x$ or $y$ is that vertex.

Suppose $x$ is the vertex; then $y$ is an end. Let $R$ be a ray of $y$ that avoids~$S$. Then the finite set $S\sub V\sm\{x\}$ separates $x$ from~$R$, a contradiction.

Suppose $y$ is the vertex. If $y\notin S$ we argue as before. Suppose that $y\in S$. Note that $y$ does not dominate~$\omega$, since $y\sim x \not\sim \omega$. But now the vertex $v\in S$ that dominates~$\omega$, if it exists, and the finitely many neighbours of ${S\sm\{v\}}$ in~$\hat C$ together separate $y$ from every ray in $x$ that avoids this finite set, a contradiction.
\end{proof}

Let us extend the notion of standard neighbourhoods from~$|G|_\Psi$ to~$\Gtilde_\Psi$. Call a neighbourhood of a point $[x]_\Psi$ of $\Gtilde_\Psi$ a \emph{standard neighbourhood} if its inverse image under $\pi_\Psi$ is a union $\bigcup_{y\in[x]_\Psi}U_y$ of standard neighbourhoods $U_y$ in $|G|_\Psi$ of the points $y\in[x]_\Psi$. Neighbourhoods in subspaces of $\Gtilde_\Psi$ that are induced by such standard neighbourhoods of~$\Gtilde_\Psi$%
   \COMMENT{}
   will likewise be called {\em standard}. All standard neighbourhoods in $\Gtilde_\Psi$ and its subspaces are open, by definition of the identification and the subspace topology.%
   \COMMENT{}

\begin{lem}\label{standardtilde}
For every point $[x]_\Psi\in\Gtilde_\Psi$ its standard neighbourhoods form a basis of open neighbourhoods in~$\Gtilde_\Psi$.
\end{lem}

\begin{proof}
Given any open neighbourhood $N$ of $[x]_\Psi$ in~$\Gtilde_\Psi$, its inverse image $W$ under $\pi_\Psi$ is open in $|G|_\Psi$%
   \COMMENT{}
   and contains every $y\in [x]_\Psi$.
By Lemma~\ref{standardnbhds}, we can find for each of these $y$ a standard neighbourhood $U_y\sub W$ of $y$ in~$|G|_\Psi$. By Lemma~\ref{stdnhd}, their union $U = \bigcup_y U_y$ is closed in~$|G|_\Psi$ under~$\sim\,$,%
   \COMMENT{}
   so $U = \pi_\Psi^{-1}(\pi_\Psi(U))$. Since $U$ is open in~$|G|_\Psi$,%
   \COMMENT{}
   this means that $\pi_\Psi(U)\sub N$ is an open neighbourhood of $[x]_\Psi$ in~$\Gtilde_\Psi$.
\end{proof}

Our next topic is to compare circles in $\Gtilde_\Psi$ with circles in~$\Gtilde$. We have already seen that circles in $\Gtilde_\Psi$ define circles in~$\Gtilde$, by composition with~$\sigma_\Psi$. The converse will generally fail: the inverse of $\sigma_\Psi$ (where it is defined) need not be continuous, so a circle in $\Gtilde$ need not induce a circle in~$\Gtilde_\Psi$ even if its points all lie in the image of~$\sigma_\Psi$. This is illustrated by the following example.

\begin{example}\label{psicontinuity}
Consider the graph of Figure~\ref{Psitrees} with $\Psi=\{\psi\}$. The closure of the double ray $D$ in $\Gtilde$ is a circle there, since in $\Gtilde$ the ends $\omega$ and~$\psi$ are identified. This circle lies in the image of~$\sigma_\Psi$, but $\sigma_\Psi^{-1}$ restricted to it fails to be continuous at the point~$\{v,\omega,\psi\}$, which $\sigma_\Psi^{-1}$ maps to the point $\{v,\psi\}$ of~$\Gtilde_\Psi$.
\end{example}

   \begin{figure}[htpb]
\centering\vskip-9pt
   	  \includegraphics[width=10cm]{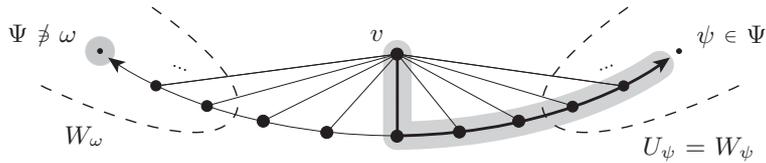}
   	  \caption{A circle in $\Gtilde$ through $p = \{v,\omega,\psi\}$ which defines for $\Psi=\{\psi\}$ a~circle in $\Gtilde_\Psi$ through~$\{v,\psi\}$.}
   \label{Circles}\vskip-6pt\vskip0pt
   \end{figure}

However, the map $\sigma_\Psi^{-1}$ in this example {\em is\/} continuous on the circle in $\Gtilde$ shown in Figure~\ref{Circles}, which `does not use' the end $\omega\in\Psi^\complement$ when it passes through the point~$\{v,\omega,\psi\}$. The fact that circles in $\Gtilde$ do induce circles in $\Gtilde_\Psi$ in such cases will be crucial to our proof of Theorem~\ref{main}:

\begin{lem}\label{circles}
~
\begin{enumerate}[\rm (i)]
   \item  Let $\rho\colon S^1\to\GPtilde$ be a circle,%
   \COMMENT{}
   with image $C$ say, and let $D$ be the set of all inner points of edges on~$C$. Then every end in the $|G|$-closure%
   \footnote{We shall freely consider $D$ as a subset of either $\GPtilde$ or~$|G|$, and similarly in~(ii).}
   of~$D$ lies in~$\Psi$.\looseness=-1
\item Let $\varphi\colon S^1\to\Gtilde$ be a circle, with image $C$ say, and let $D$ be the set of all inner points of edges on~$C$. If every end in the $|G|$-closure of~$D$ lies in~$\Psi$, then the composition $\sigma_\Psi^{-1}\circ\varphi \colon S^1\to\GPtilde$ is well defined and a circle in~${\GPtilde}$.
\end{enumerate}
\end{lem}

\begin{proof}
(i) Consider an end $\omega$ in the $|G|$-closure of~$D$. Since $|G|$ (unlike~$\Gtilde$)%
   \COMMENT{}%
   \COMMENT{}
   is first-countable, there is a sequence $(x_i)_{i\in\N}$ of points in $D$ that converges to~$\omega$ in~$|G|$.%
   \COMMENT{}
   Suppose~$\omega\in\Psi^\complement$. We show that the $x_i$ have no accumulation point on~$C$, indeed in all of~$\GPtilde$; this will contradict the fact that~$C$, being a circle, is compact and contains all the~$x_i$.

Consider a point $p\in \GPtilde$, and any representative $z\in p\sub |G|_\Psi$. As $\omega\in\Psi^\complement$ we have $\lim x_i = \omega\ne z$. Therefore $z$~has a neighbourhood $W_z$ in $|G|$ not containing any of the~$x_i$%
   \COMMENT{}
   (other than possibly $x_i=z$, which can happen only if $p=\{x_i\}$ is a singleton class).%
   \COMMENT{}
   By Lemma~\ref{standardnbhds}, the $|G|_\Psi$-neighbourhood $W_z\cap |G|_\Psi$ of $z$ contains a standard $|G|_\Psi$-neighbourhood $U_z$ of~$z$. By Lemma~\ref{standardtilde}, $\pi_\Psi\big(\bigcup_{z\in p} U_z\big)$ is a standard neighbourhood of $p$ in $\GPtilde$ that contains no $x_i$%
   \COMMENT{}
   other than possibly~$p$ itself, so $p$ is not an accumulation point of the~$x_i$.

(ii) Assume that every end in the $|G|$-closure of~$D$ lies in~$\Psi$. To show that $\sigma_\Psi^{-1}\circ\varphi$ is well defined, let us prove that $\im\varphi\sub\im{\sigma_\Psi}$. The only points of $\Gtilde$ not in the image of~${\sigma_\Psi}$ are singleton \hbox{$\sim\,$-\,classes} of $|G|$ consisting of an undominated end $\omega\notin\Psi$. By assumption and Lemma~\ref{standardnbhds}, such an end $\omega$ has a standard neighbourhood in $|G|=|G|_\Omega$%
   \COMMENT{}
   disjoint from~$D$, which $\pi$ maps to a standard neighbourhood of $\{\omega\}$ in~$\Gtilde$ disjoint from~$D$. So $\{\omega\}$ is not in the $\Gtilde$-closure of~$D$. But that closure is the entire circle~$C$, see~\cite{TST}, giving $\{\omega\}\notin\im\varphi$. This completes the proof of $\im\varphi\sub\im{\sigma_\Psi}$. As~${\sigma_\Psi}$ is injective, it follows that $\sigma_\Psi^{-1}\circ\varphi$ is well defined.

To show that $\sigma_\Psi^{-1}$ is continuous on~$C$, let a point $p\in C$ be given. Since $p$ lies in~$\im\varphi\sub\im\sigma_\Psi$, it is represented by a point $x$ in~$G\cup\Psi$; then $p=[x]$ and $\sigma_\Psi^{-1}(p) = [x]_\Psi$. By Lem\-ma~\ref{standardtilde}, it suffices to find for every standard neighbourhood $u$ of $[x]_\Psi$ in $\im (\sigma_\Psi^{-1}\!\!\restriction\! C)$ a neighbourhood $w$ of $[x]$ in~${C}$ such that $\sigma_\Psi^{-1}(w)\sub u$.

By definition, $u$~is the intersection with $\im (\sigma_\Psi^{-1}\!\!\restriction\! C)$ of a set $U\sub\Gtilde_\Psi$ whose inverse image under $\pi_\Psi$ is a union
 $$\pi_\Psi^{-1}(U) = \bigcup_{y\in[x]_\Psi}U_y$$
of standard neighbourhoods $U_y$ in $|G|_\Psi$ of the points $y\in[x]_\Psi$. Our aim is to find a similar set $W$ to define~$w$: a set $W\sub\Gtilde$ such that for $w:= W\cap{C}$ we have $\sigma_\Psi^{-1}(w)\sub u$, and such that
 \begin{equation}
\pi^{-1}(W) = \bigcup_{y\in[x]} W_y\label{Wy}
\end{equation}
where each $W_y$ is a standard neighbourhood of $y$ in~$|G|$.%
   \COMMENT{}

Let us define these~$W_y$, one for every $y\in [x]$. If $y\in G$,%
   \COMMENT{}
  then $y\in [x]_\Psi$. Hence $U_y$ is defined, and it is a standard neighbourhood of $y$ also in~$|G|$; we let ${W_y := U_y}$. If $y\in\Psi$,  then again $y\in [x]_\Psi$, and $U_y$ (exists and) has the form $\hat C_\Psi (F,y)$. We let ${W_y := \hat C_\Omega (F,y)}$ be its closure in~$|G|$; this is a standard neighbourhood of $y$ in~$|G|$. Finally, if $y\in\Psi^\complement$, then  $y\notin [x]_\Psi$ and $U_y$ is undefined. We then let $W_y$ be a standard neighbourhood of $y$ in $|G|$ that is disjoint from~$D$; this exists by assumption and Lemma~\ref{standardnbhds}. Let us call these last~$W_y$ {\em new\/}.%
   \COMMENT{}

By Lemma~\ref{stdnhd}, all these $W_y$ are closed under equivalence in $|G|\sm [y]$. Hence $\bigcup_{y\in [x]} W_y$ is closed under equivalence in~$|G|$. Its $\pi$-image $W$%
   \COMMENT{}
   therefore satisfies~\eqref{Wy} and is a standard neighbourhood of $[x]$ in~$\Gtilde$. Hence,%
   \COMMENT{}
   $w := W\cap{C}$ is a neighbourhood of~$[x]$ in~${C}$.

It remains to show that $\sigma_\Psi^{-1}$ maps every point $q\in w$ to~$u$. This is clear for $q=p=[x]$, so assume that $q\ne [x]$. By construction of~$W$ and Lemma~\ref{stdnhd}, the set $q$ lies entirely inside one of the~$W_y$.%
   \COMMENT{}
   Let us show that no such $W_y$ can be new. Since $q$ is a point in $w\sub{C}$, in which $D$ is dense~\cite{TST}, there is no neighbourhood of $q$ in $\Gtilde$ that is disjoint from~$D$. But then $q$ has an element $z$ all whose $|G|$-neighbourhoods meet~$D$. (If not, we could pick for every element of $q$ a standard $|G|$-neighbourhood disjoint from~$D$; then the union of all these would project under $\pi$ to a standard neighbourhood of $q$ in $\Gtilde$ that avoids~$D$.) As $W_y$ is a $|G|$-neighbourhood of $z\in q\sub W_y$, it thus cannot be new.%
   \COMMENT{}

We thus have $q\sub W_y$ where $W_y$ is the $|G|$-closure of~$U_y$ for some $y\in [x]_\Psi$ (or equal to~$U_y$). In particular, $W_y\sm U_y\sub\Psi^\complement$. As $q$ lies in ${C}$, in which $D$ is dense, we cannot have $q=\{\omega\}$ with $\omega\in\Psi^\complement$ (as earlier).%
   \COMMENT{}
   So either $q = \{\psi\}$ with $\psi\in\Psi$, or $q$ contains a vertex. In either case, $q\cap U_y\ne\emptyset$,%
   \COMMENT{}
   which implies that $\sigma_\Psi^{-1}(q)\in U$. As $q\in{C}$, this implies $\sigma_\Psi^{-1}(q)\in u$, as desired.
\end{proof}

	\COMMENT{}
\begin{lem}\label{closed}
Arc-components of standard subspaces of $\Gtilde_\Psi$ are closed.%
   \footnote{This refers to either the subspace or to the entire space~$\GPtilde$; the two are equivalent, since standard subspaces of~$\GPtilde$ are themselves closed in~$\GPtilde$.}
\end{lem}

\begin{proof}
Let $X$ be an arc-component of a standard subspace of $\Gtilde_\Psi$.
If $X$ is not closed, there is a point $q$ in $\Gtilde_\Psi\sm X$ such that every (standard) neighbourhood of~$q$ meets~$X$.
As in the proof of Lemma~\ref{circles},%
   \COMMENT{}
   this implies that $q$ has a representative $y\in |G|_\Psi$ such that every standard neighbourhood $U_y$ of $y$ in $|G|_\Psi$ meets~$\pi^{-1}(X)$, say in a point $x = x(U_y)$. Clearly, $y$~is an end.%
   \COMMENT{}
   Since $x\not\sim y$, we even have $[x]_\Psi\sub U_y$ by Lemma~\ref{stdnhd}.%
	\COMMENT{}
   Let $U_0\supseteq U_1\supseteq\ldots$ be a neighbourhood basis for~$y$ consisting of such standard neighbourhoods~$U_y$, and let $x_i := x(U_i)$ and $z_i := [x_i]_\Psi$ for all~$i$. Then these $x_i$ converge to $y$ in~$|G|_\Psi$, while $(z_i)_{i\in\N}$ is a sequence of points in $X$ that converges in~$\GPtilde$ to~$q=[y]_\Psi$.

For every $i\in\N\sm\{0\}$ let $A'_i$ be a $z_i$--$z_0$ arc in~$X$. Define subarcs $A_i$ of the $A'_i$ recursively, choosing as $A_i$ the initial segment of $A'_i$ from its starting point $z_i$ to its first point $a_i$ in~$\bigcup_{j<i}A_j$, where $A_0 := \{z_0\}$. (The point $a_i$ exists by the continuity of~$A'_i$, since $\bigcup_{j<i}A_j$ is closed, being a compact subspace of the Hausdorff space~$\GPtilde$.)%
   \COMMENT{}
   Note that no two $A_i$ have an edge in common.

Define an auxiliary graph $H$ with vertex set $\{A_i\mid i\in\N\}$ and edges $A_iA_j$ whenever $j$ is the smallest index less than~$i$ such that $A_i\cap A_j\ne\emptyset$. Suppose first that $H$ has a vertex $A_j$ of infinite degree. Since the arc $A_j$ is compact, it has a point $p$ every neighbourhood of which meets infinitely many~$A_i$.%
   \COMMENT{}
   By Lemma~\ref{finsepext}, there is a finite set $F$ of edges such that in $\GPtilde\sm\mathring F$ the points $p$ and~$q$ have disjoint open neighbourhoods $O_p$ and $O_q$ partitioning~$\GPtilde\sm\mathring F$. Then for infinitely many~$i$ we have both $A_i\cap O_p\ne\emptyset$ and $z_i\in O_q$. For all these~$i$ the arc~$A_i$, being connected, must have an edge in the finite set~$F$,%
   \COMMENT{}
   a contradiction.

So $H$ is locally finite. By K\"onig's infinity lemma, $H$~contains a ray $A_{i_1}A_{i_2}\dots$ such that $i_j<i_k$ whenever $j<k$.
We claim that $A:= A_{i_1}a_{i_2} A_{i_2}a_{i_3}\dots q$ is an arc in~$\GPtilde$; this will contradict our assumption that $A_{i_1}$ lies in the arc-component $X$ of $\GPtilde$ while $q$ does not.
We only have to show that $A$ is continuous in~$q$. Since every neighbourhood of $q$ in $\GPtilde$ contains the $\pi_\Psi$-image of one of our standard neighbourhoods $U_n$ of~$y$,%
   \COMMENT{}
   it suffices to show that for every such $U_n$ we have $A_i\sub \pi_\Psi(U_n)$ for all but finitely many~$i$.

Since $U_n$ is a standard neighbourhood of~$y$, there exists a set $F$ of edges such that $U_n = \hat C_\Psi (F,y)$ and $F$ is either finite or $v$-cofinite with $v\sim y$. Let $F'$ be obtained from $F$ by adding to it any other edges incident with such a vertex $v\sim y$.%
   \COMMENT{}
   Since none of the $A_i$ contains such a vertex~$v$, and distinct $A_i$ are edge-disjoint,  all but finitely many $A_i$ lie in $(\GPtilde-q)\sm\mathring F'$ and have their starting vertex $z_i = [x_i]_\Psi$ in~$\pi_\Psi(U_n)$, by the choice of~$U_n$. To complete our proof, we shall show that $\pi_\Psi (U_n\sm q)\sm\mathring F'$ and its complement in $(\GPtilde-q)\sm\mathring F'$ are two open subsets of $(\GPtilde-q)\sm\mathring F'$ partitioning it: then none of those cofinitely many $A_i$ can meet both, so they will all lie entirely in~$\pi_\Psi(U_n)$.

Since $U_n$ is a standard neighbourhood of~$y\in q$, the set $U_n\sm q$ is open in $|G|_\Psi\sm q$ and closed under equivalence, so $\pi_\Psi (U_n\sm q)$ is open in~$\GPtilde-q$ and $\pi_\Psi (U_n\sm q)\sm\mathring F'$ is open in~$(\GPtilde-q)\sm\mathring F'$. Its complement in $(\GPtilde-q)\sm\mathring F'$ is open, because it is the $\pi_\Psi$-image of the ($\sim$-closed) union of the finite set $S$ of vertices that are incident with edges in~$F$ but are not in~$U_n$,%
   \COMMENT{}
   the edges incident with them that are not in~$\mathring F'$, and the $|G|_\Psi$-closures of the components of $G-S$ not contained in~$U_n$.%
   \COMMENT{}
   The two open sets partition all of~$(\GPtilde-q)\sm\mathring F'$, because $U_n$ is itself the $|G|_\Psi$-closure of a component of $G-S$ together with the edges between $S$ and that component (which all lie in~$F$).
\end{proof}

\section{Proof of Theorem~\ref{main}}\label{proofofmain}

We can now apply the lemmas from Section~\ref{basics} to prove Theorem~2. One of these lemmas, Lemma~\ref{circles}, also implies a characterization of duality in terms of circuits and bonds. Let us include this in the statement of the theorem:

\begin{theo}
Let $G = (V, E, \Omega)$ and $G^* = (V^*\!, E, \Omega)$ be finitely separable \hbox{$2$-connected} graphs with the same edge set $E$ and the same end set~$\Omega$, in the sense of~{\rm \cite{BruhnSteinEndDuality}}. Then the following assertions are equivalent:
\begin{enumerate}[\rm (i)]
   \item $G$ and $G^*$ are duals of each other.%
   \COMMENT{}
   \item For all $\Psi\sub\Omega$ and $F\sub E$ the following holds: $F$ is the edge set of a $\Psi$-tree in~$G$ if and only if $F^\complement$ is the edge set of a $\Psi^\complement$-tree in~$G^*$.
   \item There exists a set $\Psi\sub\Omega$ such that for every $F\sub E$ the following holds: $F$ is the edge set of a $\Psi$-tree in~$G$ if and only if $F^\complement$ is the edge set of a $\Psi^\complement$-tree in~$G^*$.
   \item For all $\Psi\sub\Omega$ and $D\sub E$ the following holds: $D$~is a $\Psi$-circuit of~$G$ if and only if $D$ is a bond of~$G^*$ and every end in the closure%
   \footnote{This refers to the closure in~$|G|$ or, equivalently by~\cite{BruhnSteinEndDuality}, the closure in~$|G^*|$.}\!%
   \COMMENT{}
   of $\bigcup D$ lies in~$\Psi$.\looseness=-1
   \item There exists a set $\Psi\sub\Omega$ such that for every $D\sub E$ the following holds: $D$~is a $\Psi$-circuit of~$G$ if and only if $D$ is a bond of~$G^*$ and every end in the closure${}^{\thefootnote}\!$ of $\bigcup D$ lies in~$\Psi$.
\end{enumerate}
\end{theo}

\noindent
Remark. The fact that (i)--(iii) are symmetrical in $G$ and~$G^*$, while (iv) and~(v) are not, is immaterial and only serves to avoid clutter: as noted before, it was proved in~\cite[Theorem~3.4]{duality} that if $G^*$ is a dual of~$G$ then $G$ is a dual of~$G^*$.

\medbreak

We shall prove the implications (i)$\to$(iv)$\to$(v)$\to$(i) first, and then the implications (i)$\to$(ii)$\to$(iii)$\to$(i). The two proofs can be read independently.

(i)$\to$(iv) Assume~(i), and let $\Psi\sub\Omega$ and $D\sub E$ be given for a proof of~(iv). If $D$ is a $\Psi$-circuit of~$G$, for the circle $\rho\colon S^1\to\GPtilde$ say, it is also a circuit of~$G$%
   \COMMENT{}
   with circle $\sigma_\Psi\circ\rho\colon S^1\to\Gtilde$. By~(i), then, $D$~is a bond of~$G^*$. By Lemma~\ref{circles}\,(i), every end in the closure of $\bigcup D$ lies in~$\Psi$.

If, conversely, $D$~is a bond of~$G^*$, then $D$ is a circuit of $G$ by~(i), say with circle $\varphi\colon S^1\to\Gtilde$. If every end in the closure of $\bigcup D$ lies in~$\Psi$ then, by Lemma~\ref{circles}\,(ii), the composition $\sigma_\Psi^{-1}\circ\varphi$ is well defined and a circle in~$\GPtilde$. The edges it contains are precisely those in~$D$, so $D$ is a $\Psi$-circuit.

(iv)$\to$(v) Using the empty set for $\Psi$ in (iv) immediately yields~(v).

(v)$\to$(i) As $G$ and $G^*$ are finitely separable and 2-connected, \cite[Lemma 4.7\,(i)]{duality}%
   \COMMENT{}
   implies that $G^*$ is dual to $G$ as soon as the {\em finite\/} circuits of $G$ are precisely the finite bonds of~$G^*$. This is immediate from~(v).

\medbreak

Let us now prove the implications (i)$\to$(ii)$\to$(iii)$\to$(i). When we consider edges in $E$ topologically, we take them to include their endvertices in $\GPtilde$ or in~$\Gtilde^*_{\Psi^\complement}$, depending on the context. Thus, in (ii) and~(iii), $\bigcup F$~will be a subspace of~$\GPtilde$ while $\bigcup F^\complement$ will be a subspace of~$\Gtilde^*_{\Psi^\complement}$.

(i)$\to$(ii) We first show that (i) implies the analogue of (ii) with ordinary topological connectedness, rather than the arc-connectedness required of a $\Psi$-tree:\looseness=-1

\begin{enumerate}[($\star$)]
\item For all $F \sub E$ and $\Psi \sub \Omega$: $F$ is the edge set of a connected spanning $\Psi$-forest of~$G$ if and only if $F^\complement$ is the edge set of a connected spanning $\Psi^\complement$-forest of~$G^*$.
\end{enumerate}

For our proof of ($\star$) from~(i), let $F\sub E$ and $\Psi\sub\Omega$ be given, and assume that $F$ is the edge set of a connected spanning $\Psi$-forest $T$ of~$G$. Let $X$ be the closure in $\Gtilde^*_{\Psi^\complement}$ of $V(\Gtilde^*_{\Psi^\complement})\cup \bigcup F^\complement$. We shall prove that $X$ is a connected subspace of $\Gtilde^*_{\Psi^\complement}$ that contains no circle. Then $X$ cannot have isolated vertices, so it will be a standard subspace, and it is spanning by definition. Roughly, the idea is that $X$ should be connected because $T$ is acirclic, and acirclic because $T$ is connected.

Let us show first that $X$ contains no circle. Suppose there is a circle ${\varphi\colon S^1\to X}$, with circuit $D\sub F^\complement$ say. By Lemma~\ref{circles}\,(i) applied to $G^*$ and~$\Psi^\complement$, every end in the $|G^*|$-closure of $\bigcup D$ lies in~$\Psi^\complement$. But the ends in the $|G^*|$-closure of $\bigcup D$ are precisely those in its $|G|$-closure, by~(i). Hence we obtain:
\begin{txteq}\label{Dclosure}
The $|G|$-closure of $\bigcup D$ contains no end from~$\Psi$.
\end{txteq}

Since $D$ is also the circuit of the circle $\sigma_{\Psi^\complement}\circ\varphi\colon S^1\to \Gtilde^*$, assumption~(i) implies that $D$ is a bond in~$G$;%
   \COMMENT{}
   let $\{V_1,V_2\}$ be the corresponding partition of~$V$. Let us show the following:
\begin{txteq}\label{V1V2}
Every point $p\in \Gtilde_\Psi$ has a standard neighbourhood $N$ such that $\psi_\Psi^{-1}(N)$ contains vertices from at most one of the sets $V_1$ and~$V_2$.
\end{txteq}
Suppose $p\in \Gtilde_\Psi$ has no such neighbourhood. Then $p$ has a representative $x$ all whose standard neighbourhoods in $|G|_\Psi$ meet~$V_1$, and a representative~$y$ all whose standard neighbourhoods in $|G|_\Psi$ meet~$V_2$.%
   \COMMENT{}

If $x=y$, the point $x=y =:\psi$ is an end in~$\Psi$. Then every standard neighbourhood of $\psi$ in $|G|_\Psi$ contains a graph-theoretical path from $V_1$ to~$V_2$, and hence an edge from~$D$, because the subgraphs of $G$ underlying standard neighbourhoods in~$|G|_\Psi$ are connected and meet both $V_1$ and~$V_2$. This contradicts~\eqref{Dclosure}.

So $x\ne y$. In particular, $p$~is nontrivial, so it contains a vertex~$v$, say in~$V_1$. Then $v\ne y$, so $y =:\psi\in \Psi$. Pick a ray $R\in \psi$. Replacing $R$ with a tail of $R$ if necessary, we may assume by~\eqref{Dclosure} that $R$ has no edge in~$D$. If all the vertices of $R$ lie in~$V_1$, then every standard neighbourhood of $y=\psi$ meets both $V_1$ and~$V_2$, which contradicts~\eqref{Dclosure} as in the case of $x=y$. So $R\sub G[V_2]$. Let us show that every standard neighbourhood $\hat C_\Psi (F',\psi)$ of $\psi$ contains the inner points of an edge from~$D$, once more contrary to~\eqref{Dclosure}. 

By Lemma~\ref{standardnbhds}\,(ii), $F'$~is $v$-cofinite. Since $v\sim\psi$, there are infinitely many $v$--$R$ paths $P_0,P_1,\dots$ in~$G$ that meet pairwise only in~$v$. Since $D$ separates $v$ from~$R$, each $P_i$ contains an edge $e_i\in D$. Only finitely many of the $P_i$ contain one of the finitely many edges from $F'$ that are not incident with~$v$. All the other $P_i$ have all their points other than~$v$ in~$\hat C_\Psi (F',\psi)$, including the inner points of~$e_i$. This completes the proof of~\eqref{V1V2}.

For every point $p\in\GPtilde$ pick a standard neighbourhood $N_p$ as in~\eqref{V1V2}. Let $O_1$ be the union of those $N_p$ such that $\pi_\Psi^{-1}(N_p)$ meets~$V_1$, and $O_2$ the union of the others. Then $O_1,O_2$ are two open subsets of $\GPtilde$ covering it, and it is easy to check that $O_1\cap O_2\sub\mathring D$.%
   \COMMENT{}
So no connected subspace of $\GPtilde\sm\mathring D$ contains vertices from $V_1$ as well as from~$V_2$. But our connected spanning $\Psi$-forest $T$ is such a subspace, since its edges lie in $F\sub E\sm D$. This contradiction completes the proof that $X$ contains no circle.

For the proof of~$(\star)$ it remains to show that $X$ is connected. If not, there are open sets $O_1, O_2$ in $\Gtilde^*_{\Psi^\complement}$ that each meet~$X$ and together cover it, but intersect only outside~$X$. It is easy to check that, since $X$ contains all the vertices of~$\Gtilde^*_{\Psi^\complement}$, both $O_1$ and $O_2$ contain such a vertex but they have none in common.%
   \COMMENT{}
   For $i=1,2$, let $V_i^*$ be the set of vertices of $G^*$ representing a vertex of~$\Gtilde^*_{\Psi^\complement}$ in~$O_i$. Let $C$ be a bond contained in the cut $E(V_1^*,V_2^*)$. Note that the edges $e$ of this bond all lie in~$F$: as $e$ is connected but contained in neither~$O_i$, it cannot lie in $O_1\cup O_2 = X$.%
   \COMMENT{}
As $F$ is the edge set of a $\Psi$-forest, $C\sub F$ cannot be a $\Psi$-circuit of~$G$. By~(i), however, $C$~is a circuit of~$G$,%
   \COMMENT{}
   because it is a bond of~$G^*$.
By Lemma~\ref{circles}\,(ii), therefore, there is an end $\omega\in\Psi^\complement$ in the $|G|$-closure of~$\mathring C$; then $\omega$ also lies in the $|G^*|$-closure of~$\mathring C$.

Let us show that every standard neighbourhood $W$ of $[\omega]_{\Psi^\complement}$ in~$\Gtilde^*_{\Psi^\complement}$ contains an edge from~$C$, including its endvertices in~$\Gtilde^*_{\Psi^\complement}$. By definition, $W$~is the image under $\pi_{\Psi^\complement}$ of a subset of $|G^*|_{\Psi^\complement}$ that contains a standard neighbourhood $U$ of~$\omega$ in~$|G^*|_{\Psi^\complement}$. Since $\omega$ lies in the $|G^*|$-closure of~$\mathring C$, this $U$ either contains an edge $e\in C$ together with its endvertices in~$G^*$, or it contains one endvertex (in~$G^*$) and the interior of an edge $e\in C$ whose other endvertex dominates~$\omega$ in~$G^*$. In both cases, $e$~and its endvertices in~$\Gtilde^*_{\Psi^\complement}$ lie in~$W$.%
   \COMMENT{}

So every standard neighbourhood of $[\omega]_{\Psi^\complement}$ in~$\Gtilde^*_{\Psi^\complement}$ contains an edge from~$C$, including its endvertices in~$\Gtilde^*_{\Psi^\complement}$. In particular, it meets $X$ in both $O_1$ and~$O_2$, where this edge has its endvertices. So every neighbourhood of $[\omega]_{\Psi^\complement}$ in~$X$ meets both $O_1$ and~$O_2$. This contradicts the fact that the $O_i$ induce disjoint open subsets of~$X$ of which only one contains the point~$[\omega]_{\Psi^\complement}$. This completes the proof of~($\star$).

It remains to derive the original statement of (ii) from~($\star$). Suppose (ii)~fails, say because there is a $\Psi$-tree $T$ of~$G$, with edge set~$F$ say, such that $F^\complement$ is not the edge set of a $\Psi^\complement$-tree of~$G^*$.
By ($\star$) we know that $F^\complement$ is the edge set of a connected spanning $\Psi^\complement$-forest $X$%
   \COMMENT{}
   in~$G^*$, which we now want to show is even arc-connected. Suppose it is not. Since the arc-components of $X$ are closed (Lemma~\ref{closed}), no arc-component of $X$ contains all its vertices.%
	\COMMENT{}
Vertices in different arc-components are joined by a finite path in~$G^*$, which contains an edge $e$ whose endvertices lie in different arc-components of~$X$. Then $X\cup e$ still contains no circle, so $F^\complement\cup\{e\}$ too is the edge set of a connected spanning $\Psi^\complement$-forest of~$G^*$.%
   \COMMENT{}
   Thus, by~($\star$), $F\sm\{e\}$ is the edge set of a connected spanning $\Psi$-forest of~$G$. This can only be $T\sm\mathring e$, so $T\sm\mathring e$ has precisely two path components $D_1$ and~$D_2$ but is still connected. Then $D_1$~and $D_2$ cannot both be open, or equivalently, cannot both be closed. This contradicts Lemma~\ref{closed}.

(ii)$\to$(iii) Using the empty set for $\Psi$ in (ii) immediately yields~(iii).

(iii)$\to$(i) As $G$ and $G^*$ are finitely separable and 2-connected, it suffices by \cite[Lemma 4.7\,(i)]{duality}%
   \COMMENT{}
   to show that $G^*$ is a finitary dual of $G$, i.e., that the finite circuits of $G$ are precisely the finite bonds of~$G^*$. By Lemma~\ref{matroidlemma}\,(ii), a finite set $F$ of edges is a bond of $G^*$ if and only if it meets the edge set of every $\Psi^\complement$-tree of~$G^*$ and is minimal with this property. By~(iii), this is the case if and only if $F$ is not contained in%
   \COMMENT{}
   the edge set of any $\Psi$-tree of~$G$, and is minimal with this property. By Lemma~\ref{matroidlemma}\,(i), this is the case if and only if $F$ is a circuit of~$G$.\qed

\section*{Acknowledgement}
The ideas that led to the formulation of Theorem~2 were developed jointly with Henning Bruhn. We benefited greatly from his insights at this stage.

\bibliographystyle{amsplain}
\bibliography{collective}
\small
\vskip2mm plus 1fill
Version 7 June, 2011

\end{document}